\documentclass[leqno,a4paper,oneside]{amsart}
\usepackage{amsmath,amssymb,amscd,mathrsfs,mathtools}
\usepackage{graphicx,comment,multirow,color,caption}
\usepackage[%
    font={small,sf},
    labelfont=normalfont,
    format=hang,    
    format=plain,
    margin=0pt,
    width=0.8\textwidth,
]{caption}
\usepackage[list=true]{subcaption}
\usepackage{wrapfig}
\usepackage{mathptmx,courier}
\usepackage[scaled=1.0]{helvet}
\title[Enumeration of doubly symmetric diagrams]{The enumeration of doubly symmetric diagrams for strongly positive amphicheiral knots}

\author{Christoph Lamm}

\theoremstyle{plain}
  \newtheorem{theorem}{Theorem}[section]

\theoremstyle{definition}
  \newtheorem{definition}[theorem]{Definition}
  \newtheorem{question}[theorem]{Question}

\begin{document} %%%%%%%%%%%%%%%%%%%%%%%%%%%%%%%%%%%%%%%%%%%%%%%%%%%%%%%%%%%%
%%%%%%%%%%%%%%%%%%%%%%%%%%%%%%%%%%%%%%%%%%%%%%%%%%%%%%%%%%%%%%%%%%%%%%%%%%%%%

\begin{abstract}
This is the second part of the article on doubly symmetric diagrams and strongly positive amphicheiral knots.
We develop an enumeration strategy for prime knots given by doubly symmetric diagrams and determine
all cases up to 18 crossings in the doubly symmetric diagram. A digression covers the origin of Gauss words 
for long curves and explains how Gauss marked non-realizable words.
\end{abstract}

\keywords{strongly positive amphicheiral knots, symmetric unions of knots, Gauss words}
\subjclass[2020]{57K10}

%%%%%%%%%%%%%%%%%%%%%%%%%%%%%%%%%%%%%%%%%%%%%%%%%%%%%%%

% captions seem to be too near to the tables
\captionsetup{belowskip=8pt,aboveskip=8pt}

% use this if 'oneside' is active
\reversemarginpar

%%%%%%%%%%%%%%%%%%%%%%%%%%%%%%%%%%%%%%%%%%%%%%%%%%%%%%%

\maketitle

\section{Introduction} \label{sec:Introduction}

The present article is the second part of a project concerning doubly symmetric knots. 
The first part is the article \cite{Lamm2023}. In the first part we used the template notation
for doubly symmetric diagrams in a practical way: The aim was to list the strongly positive 
amphicheiral knots we had found to be doubly symmetric and to analyze their properties.

We now provide the theoretical foundation for the relationship between templates and doubly 
symmetric diagrams. A systematic enumeration of templates enables us to answer questions 
about the `doubly symmetric' crossing number.

\begin{definition}
For a knot $K$ we denote, as usual, by $c(K)$ the minimal crossing number of a diagram of $K$.
If $K$ has a doubly symmetric diagram, we denote by $c_{ds}(K)$ the minimal crossing number of 
all doubly symmetric diagrams of $K$. This invariant is called the doubly symmetric crossing number.
\end{definition}

This article answers the following question:

\begin{question} \label{q1}
Which prime knots have $c_{ds} \le 18$?
\end{question}

This question will be answered by a systematic enumeration of all templates up to 3 crossings which are
not reducible to less complicated ones. Initially we thought that there should be a doubly symmetric 
diagram with 12 crossings which represents a prime knot. The enumeration answers the following question negatively:

\begin{question} \label{q2}
Does a non-trivial prime knot exist with $c_{ds} \le 12$?
\end{question}

Starting with 14 crossings, there are prime doubly symmetric knots. 
For instance, the prime knots with $c_{ds} = 14$ are $10_{99}$, $10_{123}$, 12n706 and 14n9732. 
The case 14n9732 also shows that the minimal crossing numbers $c(K)$ and $c_{ds}(K)$ can be equal.

\subsection{Doubly symmetric diagrams and their templates}
We assume that the definitions of doubly symmetric diagrams and their templates are known from \cite{Lamm2023},
and, in addition, use the following terms: An integer marking, symbolizing a twist box on an axis,
consists of a \textsl{twist marker} (the dotted line), and an integer, the \textsl{twist number}.
Because the enumeration approach uses templates, the first step is the following theorem:

\begin{theorem} \label{diagram_representation}
There is a bijection between the set of all doubly symmetric diagrams and the set of
templates modulo addition of adjacent twist numbers of the same sign.
\end{theorem} 

\begin{figure}[hbtp]
\centering
\includegraphics[scale=0.8]{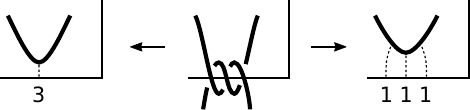}
\caption{Illustration of the ambiguity in the template construction: A twist box with $|n_i| > 1$
can be subdivided in several ways. The example shows one twist marker with twist number 3 (left), and three 
adjacent twist markers with twist numbers 1 (right). The same observation holds for crossings on the y-axis.}
\label{diagram_to_template}
\end{figure}

\begin{figure}[hbtp]
\centering
\includegraphics[scale=0.8]{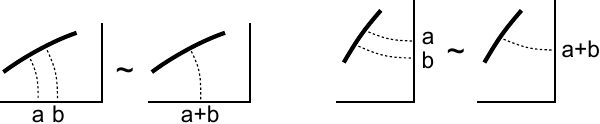}
\caption{Addition rules for adjacent twist markers, if the twist numbers $a$ and $b$ have the same sign.} 
\label{twist_markers_addition}
\end{figure}

\begin{proof}
Besides planar isotopies we do not allow any other equivalence relation for doubly symmetric diagrams.
On the other hand, the construction step from doubly symmetric diagram to template is not unique:
Crossings on one of the axes are represented as twist numbers $n_i$. If $|n_i| > 1$, the complete twist
tangle with $n_i$ half-twists can be read off in one go, or subdivided into several sub-tangles, 
see Figure \ref{diagram_to_template}. Therefore, we allow the addition of adjacent twist numbers 
of the same sign, illustrated in Figure \ref{twist_markers_addition}.
\end{proof}

\vspace{-0.5cm}

Note, that adjacent positive and negative twist numbers cannot be added, because that would 
correspond to the application of Reidemeister 2 moves in the doubly symmetric diagram.

\begin{figure}[hbtp]
\centering
\includegraphics[scale=0.8]{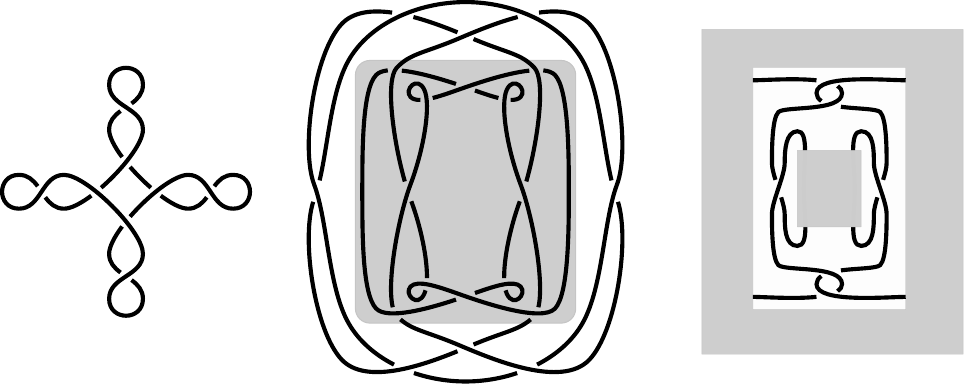}
\caption{Doubly symmetric diagrams with removable curls and crossings.}
\label{diagrams_reducible}
\end{figure}

\begin{figure}[hbtp]
\centering
\includegraphics[scale=0.75]{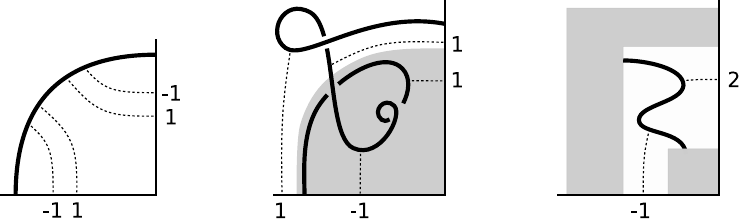}
\caption{The templates for the diagrams in Figure \ref{diagrams_reducible}. 
They illustrate cases with twist markers connected to central arcs, see Definition \ref{arc_types}.}
\label{diagrams_template_reducible}
\end{figure}

Figures \ref{diagrams_reducible} and \ref{diagrams_template_reducible} show three diagrams and 
their templates. The left example illustrates that the twists numbers $1$ and $-1$ on one of the
axes cannot be added. The example in the middle shows a template with an isolated Reidemeister 1
curl, occurring four times in the doubly symmetric diagram. This example and the one on the right
also allow the removal of crossings by rotating the central part.

For our goal of enumerating minimal diagrams for prime doubly symmetric knots, we do 
not include cases as given in these examples. This is explained in the next section.

\subsection{A first glance at the enumeration strategy}

\begin{figure}[hbtp]
\centering
\includegraphics[scale=0.8]{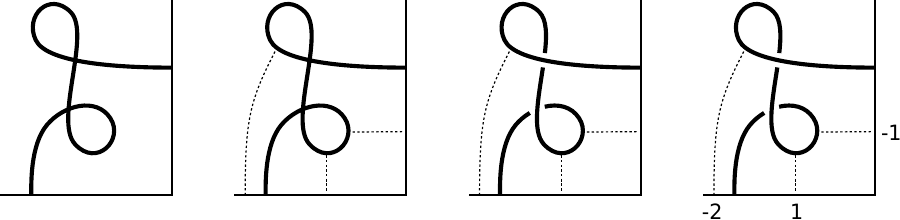}
\caption{Template types, from left to right: long planar curve, parterre template,
knot diagram templates without and with twist numbers.}
\label{templates_four_types}
\end{figure}

Our systematic enumeration proceeds in several steps, starting with the simplest template parts
and adding complexity in each step:

\begin{enumerate}
\item The basic ingredient of a template is a \textit{long planar curve} in the second quadrant of a coordinate system. 
This starts at a point on the x-axis and ends at a point on the y-axis. The start and end points are required to be different from the origin.
The number of crossings in the long curve is denoted by $n$.
An example with $n=2$ is shown in the first illustration in Figure \ref{templates_four_types}.
\item Adding twist markers (without twist numbers) yields a \textit{parterre template} (second illustration in Figure \ref{templates_four_types}).
\item With an additional choice of over- and under-crossings for the long curve we obtain a
\textit{knot diagram template} (not yet with specified twist numbers, third illustration in Figure \ref{templates_four_types}).
We also call this a \textit{parterre template with crossing variation}.
\item When the twist numbers are specified we have a \textit{knot diagram template with twist numbers}
(fourth illustration in Figure \ref{templates_four_types}). This represents a doubly symmetric knot diagram. 
\end{enumerate}

With this approach, we may sort out those variations in an early step of the enumeration, which do not contribute 
to our desired list of prime knot diagrams with a minimal number of crossings. For instance, the template in the middle of Figure 
\ref{diagrams_template_reducible} has a Reidemeister 1 curl in the interior of a diagram region and the curl cannot be connected by twist markers 
to one of the axes. It will stay therefore in any later step in this removable form, and we may already in the `long curve' step conclude that 
this curve cannot yield a doubly symmetric diagram for a prime knot with minimal $c_{ds}$.
Similarly, we sort out diagrams which allow simplifications with Reidemeister 2 moves on one of the axes 
(see the template in the left of Figure \ref{diagrams_template_reducible}).

\enlargethispage{1cm}

\subsection{Plan of the article}
The plan of the article is: 
\begin{itemize}
\item[a)] We present diagram moves and filter rules for long curves and parterre templates. A long curve or parterre 
  template which cannot yield a doubly symmetric diagram for a prime knot with minimal $c_{ds}$ is discarded. 
	The list of long curves is generated and the filters for parterre diagrams are applied for $n = 0$, 1 and 2 self-intersections.
\item[b)]	We then dive into a historical digression on long curves, giving some details on the work of 
  Carl Friedrich Gauss (in 1825), Gyula Sz\"{o}kefalvi Nagy (in 1925) and Vladimir Arnold (in 1994).
\item[c)] The list of long curves is generated for curves with $n = 3$ and we apply the filters. 
  For $n \le 3$ this results in 9 parterre templates, see Figure \ref{ParterreTemplatesNineCases}.
\item[d)] Choices of over- and under-crossings are analyzed, and for the final list of knot diagram templates we use Knotscape to identify the knots.
\item[e)] Analysis of diagrams for $n = 4$. For up to 3 self-intersections the number of cases in the sections c) and d) is reasonably small. 
  For 4 self-intersections it is much larger and we do not give a complete analysis of all cases.
\end{itemize}

\begin{figure}[hbtp]
\centering
\includegraphics[scale=0.9]{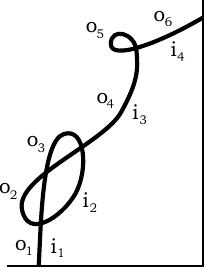} \hspace{2cm}
\includegraphics[scale=0.9]{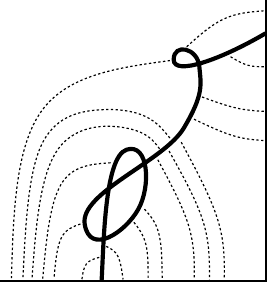}
\caption{Left: Arcs in a template are marked with $o_1$, $o_2$, \ldots (adjacent to the outer region), 
$i_1$, $i_2$, \ldots (adjacent to the inner region), and without marking (not adjacent to the outer or inner region). 
In this example, there are three arcs which bound both, the outer and inner, regions. These are called central arcs. 
Right: Example for a parterre template with twist markers.} 
\label{template_regions}
\end{figure}

\begin{definition} \label{arc_types}
For a template containing a long curve, we use the following notation for the arcs (curve parts between crossings, or a crossing 
and x- or y-axis): Outer arcs are adjacent to the outer region, and inner arcs are adjacent to the inner region. Arcs which are 
adjacent to both regions are also called central arcs. Outer and inner arcs are numbered with an index as shown in Figure \ref{template_regions}, 
starting from the x-axis and ending at the y-axis (we need this index for the definition of `separators' in a later section). 
Templates may also contain arcs which are neither outer, nor inner arcs.
\end{definition}

\section{Diagram moves, additional principles and filter rules}
We introduce three techniques that help reduce the number of enumerated cases.

\subsection{Diagram moves}
\begin{figure}[hbtp]
\centering
\includegraphics[scale=0.82]{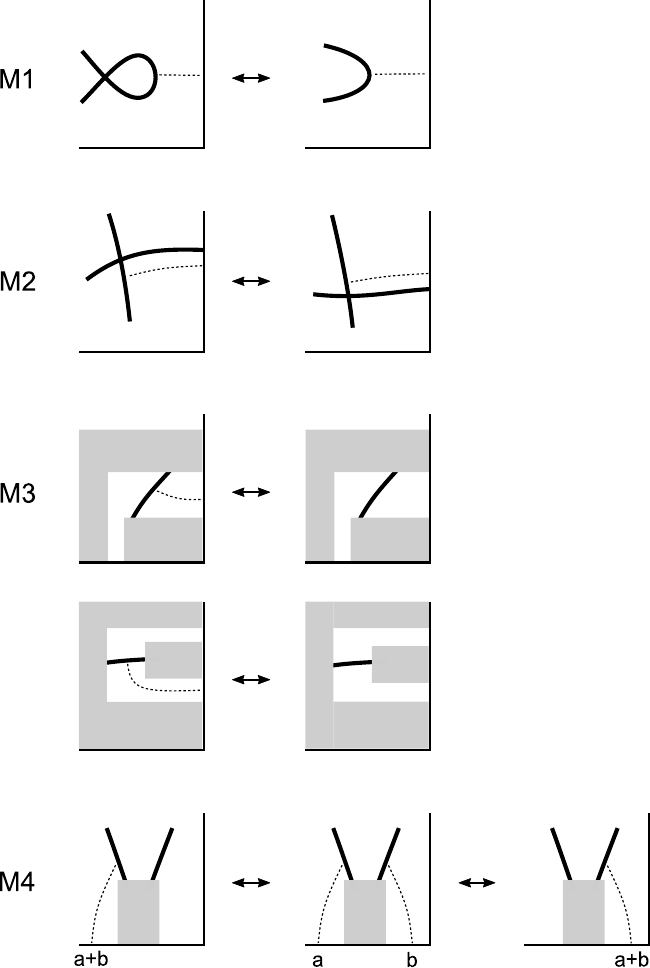}
\caption{The diagram moves} 
\label{diagram_moves}
\end{figure}

Similar to Reidemeister moves, we identify four moves for templates, preserving the knot type. 
These diagram moves are useful in our description of filters. The diagram moves M1 and M2 can be applied with respect 
to the y-axis only; the moves M3 and M4 use rotations of diagram parts and they are possible with respect to both axes
(the illustration in Figure \ref{diagram_moves} contains only one of the two versions in each case).

We first explain diagram moves, and then additional principles and filter rules.
The diagram moves are shown in Figure \ref{diagram_moves}, and are defined as follows:

\begin{itemize}
	\item M1, `Twist box flip':
				If a Reidemeister 1 curl is connected to the y-axis by a twist marker, but has no connection 
	      to the x-axis, then a flip of the twist box is possible, removing the two curls on each side of the y-axis.
				The same flip can be applied to the lower part of the diagram (quadrants III and IV), of course.
	\item M2, `One or several Reidemeister 3 moves':
	      An arc, ending at the y-axis, is moved over (or under) a twist box, 
	      using Reidemeister 3 moves. This can also be viewed as moving the twist marker under (or over) a crossing. 				
	\item M3, `Central arc flip':
	      A twist marker ending at a central arc (bounding both the inner and the outer regions) represents
	      twists on the x- or y-axis which can be removed by one or several flips. In the examples in Figure \ref{diagrams_reducible} 
				(middle and the right), the gray box can be rotated, removing crossings on the axes.		
	\item M4, `Addition of twist markers':
	      In this situation, one twist marker can be omitted, because flipping the central diagram part results in the addition of the indicated twist numbers.
\end{itemize}

\bigskip

\subsection{Additional principles}
We have the following three additional principles which reduce the number of diagrams.
These are applied---in ascending template complexity---for long curves (P1), parterre templates (P2) and knot diagram templates (P3).

\begin{itemize}
	\item P1, `Switch inner and outer regions': \\ 
	      Two diagrams which are related by this operation yield the same knot. 
	      Therefore, we need to consider only one of them. This leads to a reduction by a factor of 2 in the set of long curves.
				An example, on the detail level of knot diagrams, is shown in Figure \ref{crossings_12_composite}.
	\item P2, `Zero, one or two twist markers for each arc': \\ 
	      To produce minimal diagrams, it is sufficient to have one twist marker for each outer or inner (non-central) arc, because adjacent 
				twist numbers of different signs can be simplified by Reidemeister 2 moves. This general rule is modified as follows, however: Arcs 
				with a `separator' may have one twist marker ending at the x-axis, and one ending at the y-axis, see Definition \ref{separator}.
				We do not connect twist markers to central arcs, because these can be removed by central arc flips, M3.				
	\item P3, `Crossing choice': \\ 
	      For knot diagram templates, we choose over- or under-crossings for the long curve part. 
	      Choices which are mirror images of each other lead to the same knots, since the knots we generate are amphicheiral. 
				Therefore, we have a factor 2 reduction as well.
\end{itemize}

\bigskip				
\begin{figure}[hbtp]
\centering
\includegraphics[scale=0.6]{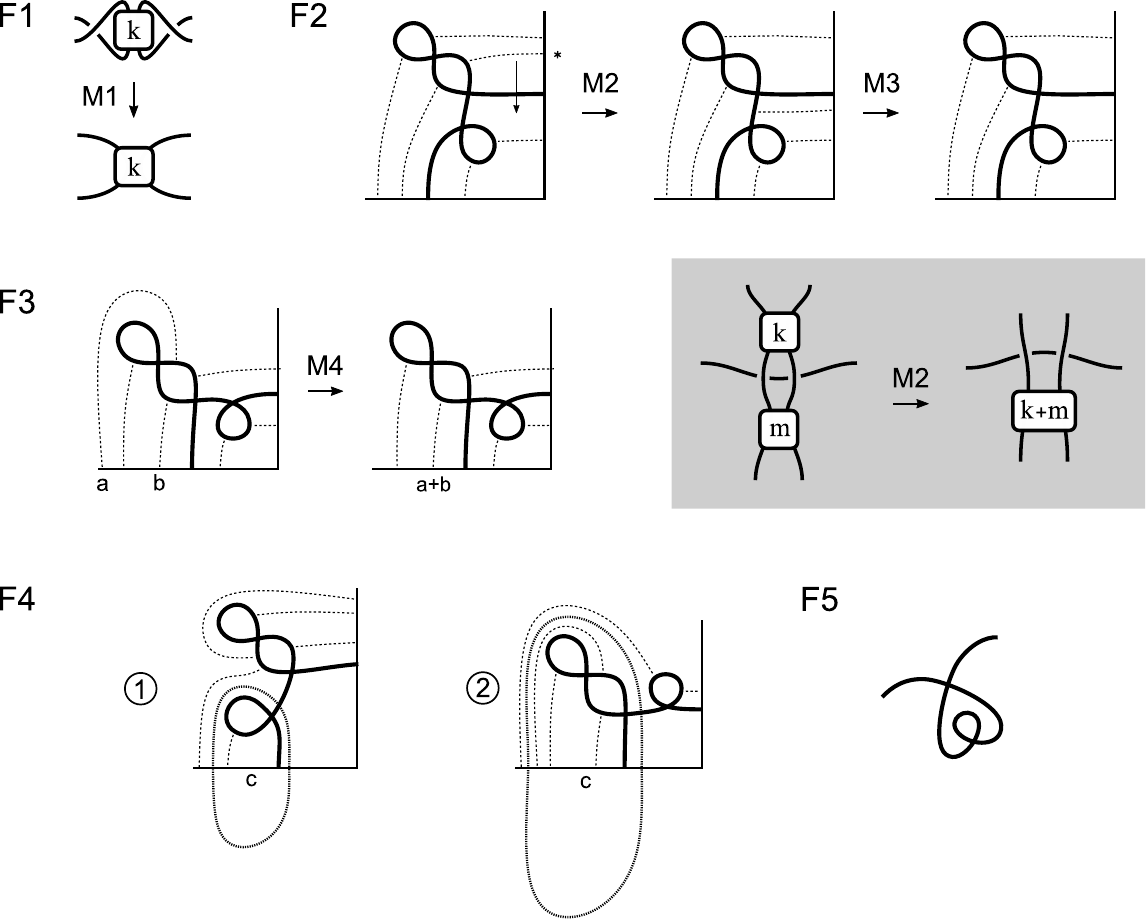}
\caption{Illustration for the filter rules. The filters F1, F4 and F5 sort out complete templates. Filters F2 and F3 
allow the omission of a twist marker. The gray box illustrates the addition of twist numbers following an M2 move.} 
\label{filter_examples}
\end{figure}

\subsection{Filter rules}
The term filter rules denotes diagram conditions which allow us to drop those templates in 
an early step of the enumeration, which cannot contribute to our desired list of knot diagrams. 
We define five filter rules; they are illustrated in Figure \ref{filter_examples}.

\begin{itemize}
  \item F1 applies, if an M1 move is possible, which reduces the number of crossings. In this case, the template is omitted.
	\item F2 is based on a combination of the moves M2 and M3.
	      A twist marker is modified by an M2 move, so that it ends at a central arc. It then allows an M3 move and can be omitted.				
	\item F3 uses the move M4. If this rule is applied, one of the two twist markers, with labels a and b, can be omitted.
	\item F4 covers composite diagrams and does not use one of the moves. If a template has composite diagram parts---contradicting our 
	      enumeration goal of generating prime knots---we omit the template. The rule applies, if the template contains a half-circle, 
				ending at the x- or y-axis, which intersects the long curve once and does not meet one of the twist markers. 
				The letter `c' indicates this composite case (and the circle).
  \item F5: If a template contains a Reidemeister 1 curl in the interior of a diagram region, we omit the template
			  because the curl cannot be connected by twist markers to one of the axes.
\end{itemize}
				
\def\nineCases{0.6}

\begin{figure}[hbtp]
\centering
\subcaptionbox{X1}{\includegraphics[scale=\nineCases]{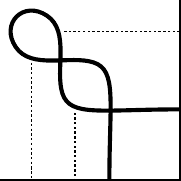}} \quad
\subcaptionbox{X2}{\includegraphics[scale=\nineCases]{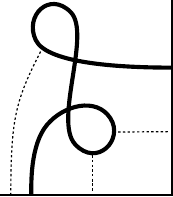}} \quad
\subcaptionbox{H1}{\includegraphics[scale=\nineCases]{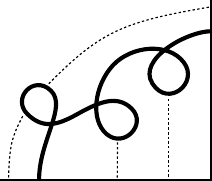}} \quad
\subcaptionbox{H2}{\includegraphics[scale=\nineCases]{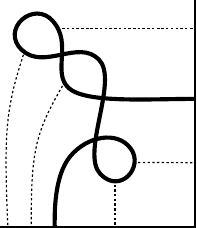}} \quad
\subcaptionbox{H3}{\includegraphics[scale=\nineCases]{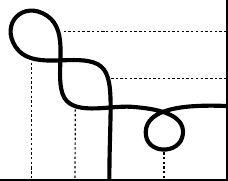}} \\
\subcaptionbox{H4}{\includegraphics[scale=\nineCases]{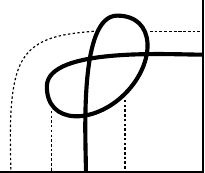}} \quad
\subcaptionbox{H5}{\includegraphics[scale=\nineCases]{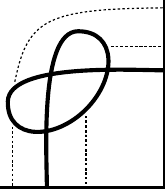}} \quad
\subcaptionbox{H6}{\includegraphics[scale=\nineCases]{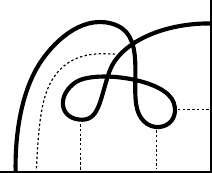}} \quad
\subcaptionbox{H7}{\includegraphics[scale=\nineCases]{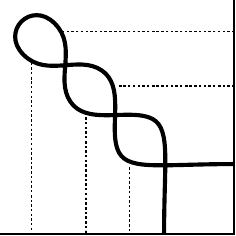}}
\caption{The nine cases for parterre templates with $n \le 3$.} 
\label{ParterreTemplatesNineCases}
\end{figure}

When the filters and the additional principles are applied, we find the parterre templates for up to three crossings shown in Figure \ref{ParterreTemplatesNineCases}. 
In Sections \ref{result_n_2} and \ref{result_n_3}, we prove that these generate all non-trivial prime knots given by doubly symmetric diagrams with $n \le 3$.

We will also give the list of generating templates for $n=4$ (which is much longer).
Before our digression on long curves (Gauss, Nagy, Arnold) we apply the filters and the additional principles for the long curves with $n \le 2$.

\subsection{The case $n=0$}
The long curve consists of a single arc without self-intersections, as in the left of Figure \ref{diagrams_template_reducible}.
All twist markers may be removed by central arc flips, M3. Therefore, in the case $n=0$ we obtain the trivial knot only.

\bigskip

\begin{figure}[hbtp]
\centering
\includegraphics[scale=0.6]{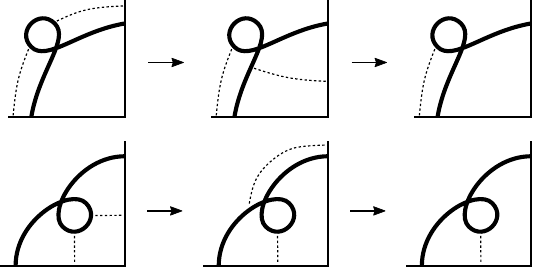}
\caption{There are two cases with $n=1$. The F2 filter is applied. The second column illustrates the M2 move 
and the third column the M3 move, which allows us to omit the twist markers connected to the y-axis.} 
\label{template_n_1}
\end{figure}

%\clearpage
%\newpage

\subsection{The case $n=1$}
There are two long curve possibilities, see the upper and lower rows in Figure \ref{template_n_1}, where the filter application is described
(we did not apply the principle P1, `Switch inner and outer regions', here and show both possibilities).
We obtain two diagrams without twist markers connected to the y-axis, meaning that the symmetric unions are composite of the form $K \sharp -K$. 
Alternatively, we could arrive at this assessment by applying the filter F4 to the diagrams in the third column.

\subsection{The case $n=2$} \label{Case_n_two_firstPart}
\begin{figure}[hbtp]
\centering
\includegraphics[scale=0.7]{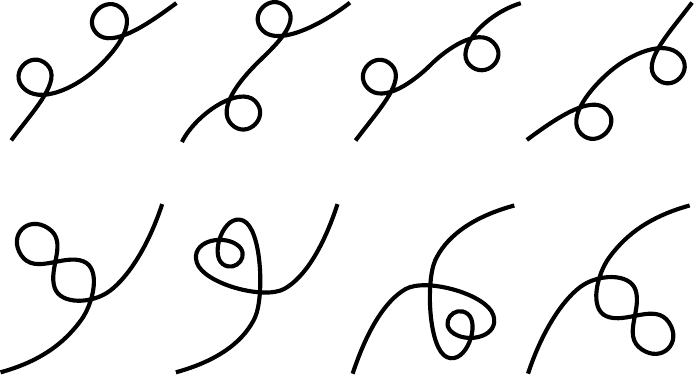}
\caption{The eight cases for long curves with $n=2$.} 
\label{template_n_2}
\end{figure}

The eight long curves with $n=2$ are shown in Figure \ref{template_n_2}. We omit the coordinate system.
The curves are oriented from lower left side to upper right side.

\begin{wrapfigure}{o}[2.5cm]{0.5\textwidth}
\flushright
\vspace{-0.7cm}
\includegraphics[width=0.5\textwidth]{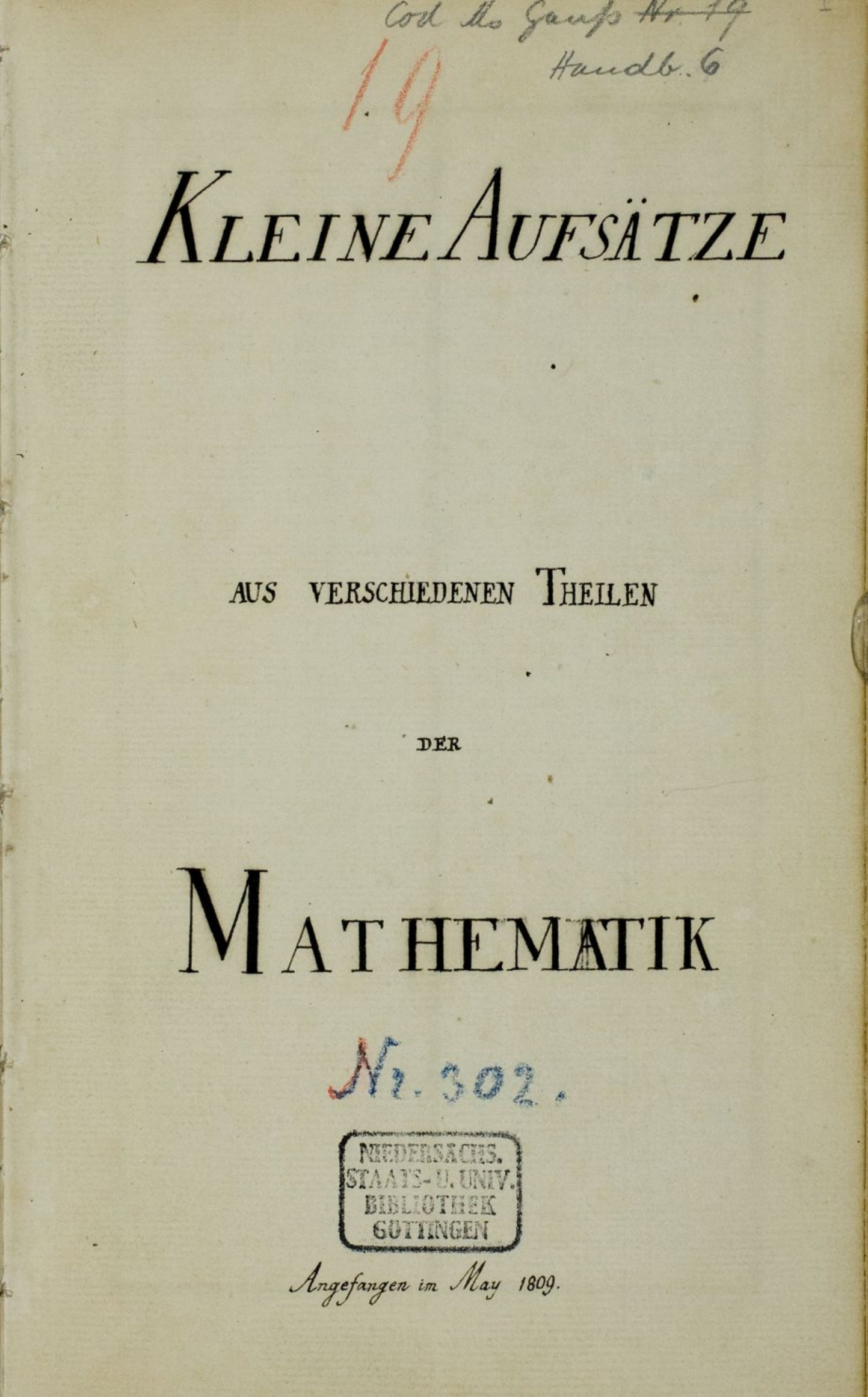}
\caption{Gauss: Title page of handbook 6} 
\label{Handbuch_6_Titelseite}
\end{wrapfigure}

The first row contains the long curves with Gauss word aabb and the second row those with abba.
The principle P1, ’Switch inner and outer regions’, can be applied for each row, and we also filter 
out the curves with a Reidemeister 1 curl in the interior of a diagram region, F5. 
Therefore, for the first row we may choose the first two curves and for the second row the first curve.

The analysis of the filter steps containing twist markers is more complicated than for $n \le 1$ and will be done in Section \ref{result_n_2}.
The result is: The two templates X1 and X2 in Figure \ref{ParterreTemplatesNineCases} generate all non-trivial prime 
knots given by doubly symmetric diagrams with $n = 2$ (and hence with $n \le 2$). 

\section{Historical remarks on long curves}
\begin{figure}[hbtp]
\centering
\includegraphics[scale=0.45]{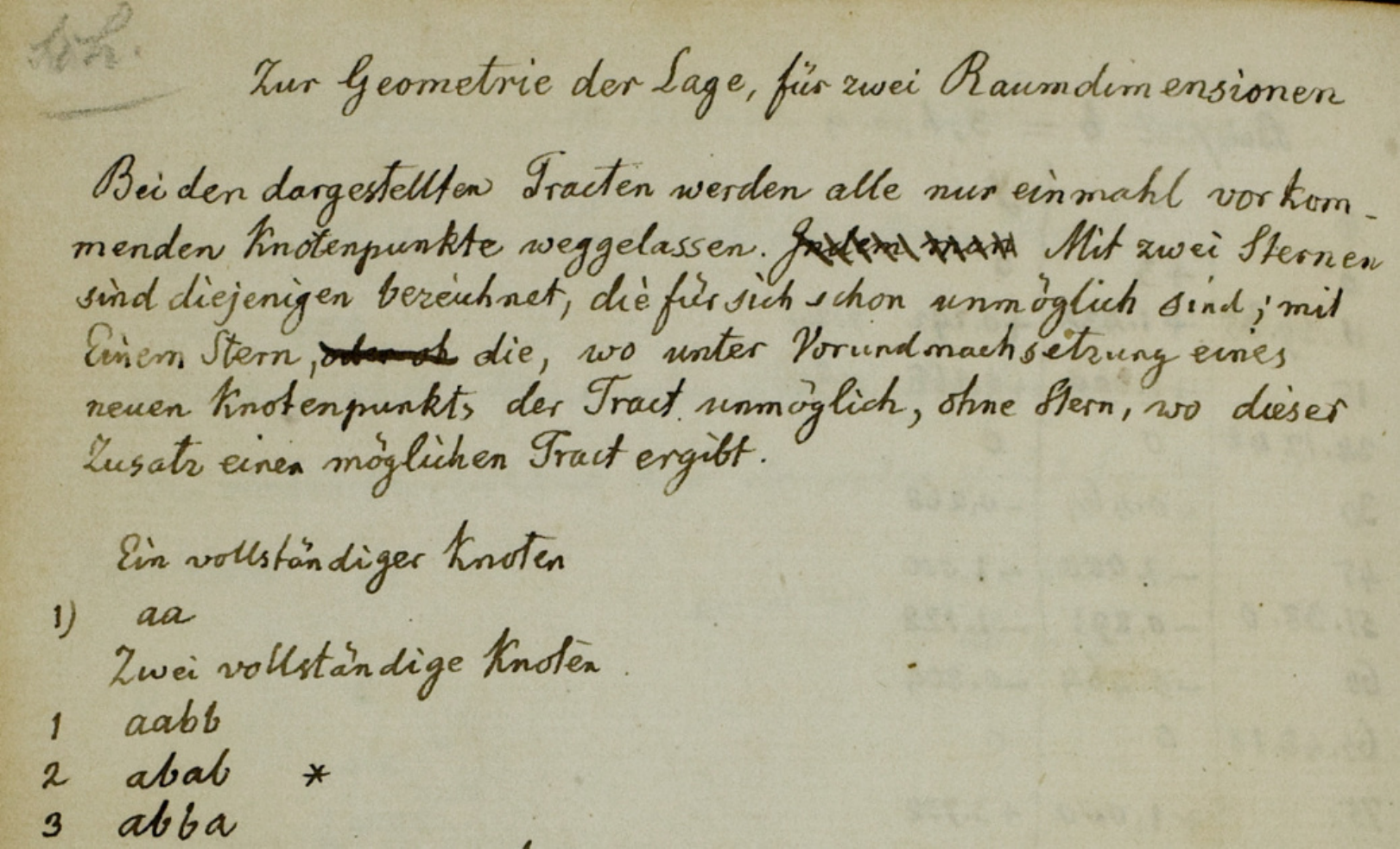}
\caption{The beginning of the section on `Tractfiguren'} 
\label{gauss_tracts_words}
\end{figure}

Carl Friedrich Gauss kept notebooks with unpublished results. His work on knot projections is part of a volume with
the title page shown in Figure \ref{gauss_tracts_words}, \cite{Gauss_Handbuch_6}.

This volume is also called `Handbuch 6', but the title, probably handdrawn by Gauss, is `Kleine Aufs\"{a}tze aus
verschiedenen Theilen der Mathematik -- Angefangen im May 1809' (Short essays on several parts of mathematics --
begun in May 1809). The formal title page indicates that this handbook contains final drafts. 

Compare this with handbook 7 which contains, for instance, raw data on geodesic topics and also the well-known braid 
sketch, \cite{Gauss_Handbuch_7}. That handbook does not have a title page.

Some of the material in the notebooks was chosen posthumously to be included in his collected works.

For instance, in volume 8 of his `Werke', \cite{Gauss_VIII}, we find the section `Geometria Situs. Nachtr\"{a}ge zu Band IV' (Addendum to Volume IV). 
Its second part is about the enumeration of knot projections and long curves and is called  `Zur Geometrie der Lage, f\"{u}r zwei Raumdimensionen'. 
Gauss called the two-dimensional curves `Tracte' or `Tractfiguren', see Figure \ref{Handbuch_6_Titelseite}.
Compare this also with the printed pages in \cite{Gauss_VIII}.

The case $n=1$ is denoted by `Ein vollst\"{a}ndiger Knoten' (one complete knot) and corresponds to
the long curves we have shown in Figure \ref{template_n_1}.

\begin{figure}[hbtp]
\centering
\includegraphics[scale=0.8]{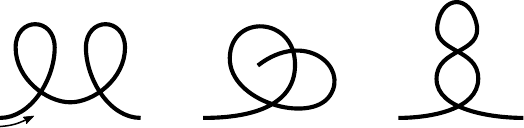}
\caption{Long curve examples for the Gauss words aabb, abab and abba. Gauss marked the second case with one star.} 
\label{gauss_words}
\end{figure}

The case $n=2$, `Zwei vollst\"{a}ndige Knoten' (two complete knots), has the three entries we illustrate in Figure \ref{gauss_words}.
The first and third curve already appeared in Figure \ref{template_n_2} together with variations.
The second one is not realizable as a closed curve. Gauss used a notation with one or two stars to mark these cases.
Non-realizable Gauss words have prompted the study of virtual knots, \cite{Kauffman}.

\subsection{Star markings} 

\begin{wrapfigure}{o}[2cm]{0.37\textwidth}
\flushright
\vspace{-0.8cm}
\includegraphics[width=0.32\textwidth]{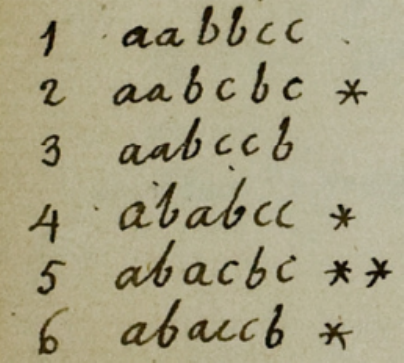}
\label{gauss_n_3_part}
\end{wrapfigure}

We discuss the meaning of the markings with one or two stars. Gauss words with two stars occur in the list with 
15 possibilities for $n=3$: It contains 6 words without a star, 7 words with one star, and 2 words with two stars,
see the examples on the right.

Gauss first explains the use of two stars (occurring for abacbc and abcbca): 
`Mit zwei Sternen sind diejenigen bezeichnet, die für sich schon unmöglich sind.' (A marking with two stars is
used for those words that are impossible per se.);

Then, he explains the use of one star: `Mit einem Stern die, wo unter Vor- und Nachsetzung eines neuen Knotenpunkts 
der Tract unmöglich [ist].' (Words are marked with one star if they are impossible when one label is added in the 
beginning and one is added at the end.); 

And finally those without a star: `Ohne Stern, wo dieser Zusatz einen möglichen Tract ergibt.' 
(Without a star if this adding of labels yields a possible word.) 

The use of the terms `Knoten' and `Knotenpunkt' is not very well defined here and we translate it with either label 
(meaning one of the letters used for the crossings), or crossing. 

The procedure of adding a crossing replaces a word $w$ by a new word $xwx$ with a new label $x$.
If the labels are traversed in alphabetical order, then we would need to shift every label by one ($a \mapsto b$,
$b \mapsto c$, $\ldots$) and then replace the shifted word $w'$ by $aw'a$.

\begin{wrapfigure}{o}[2.2cm]{0.12\textwidth}
\flushright
\vspace{-0.8cm}
\includegraphics[width=0.1\textwidth]{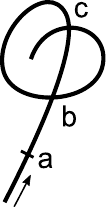}
\label{gauss_word_abcbca}
\end{wrapfigure}

\begin{figure}[hbtp]
\centering
\includegraphics[scale=0.8]{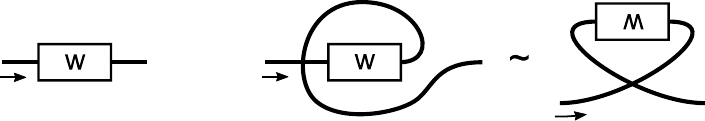}
\caption{Explanation for `Vor- und Nachsetzung eines neuen Knotenpunkts'} 
\label{gauss_word_enclosure}
\end{figure}

The case that the result is a possible word (Gauss used no star marking) is explained in Figure \ref{gauss_word_enclosure}.
In the case that the word $w$ is already marked with one star, the end of the long curve is trapped in an inner region
and therefore the adding of a crossing is not possible (abab in Figure \ref{gauss_words}).
The word $abcbca$, mentioned above, is of this type ($w'=bcbc$ and label $a$ added, see the illustration on the right). 
To better understand the star markings, the reader is invited to sketch the other example with two stars and $n=3$, $abacbc$.

\subsection{Four crossings} 
For four crossings there are already 105 words. These are listed by Gauss in 
alphabetical order and marked with one or two stars, as described in the last section.

\begin{figure}[hbtp]
\centering
\includegraphics[scale=0.51]{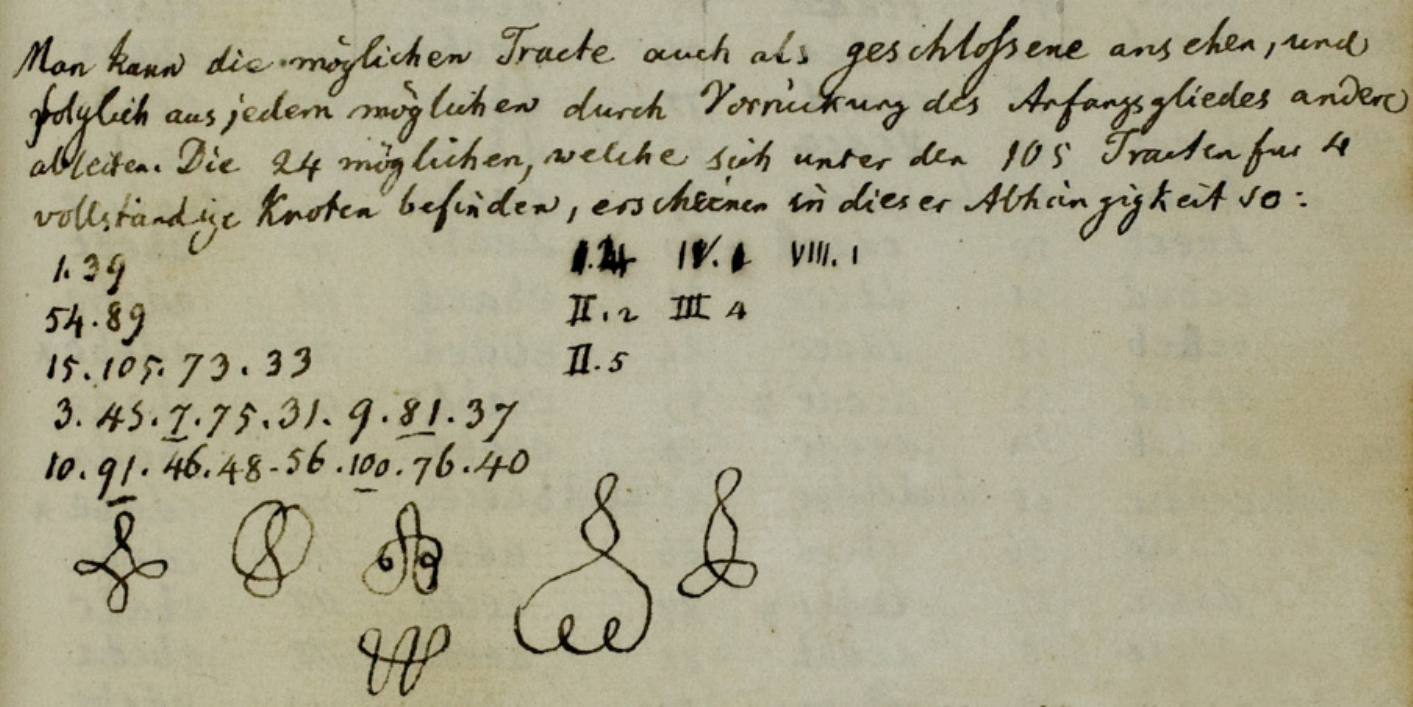}
\caption{The section in which Gauss explains the correspondence between long curves and closed curves (for $n=4$)} 
\label{Tractfiguren_handschriftlich}
\end{figure}

He finds 24 possible words and goes on with `Man kann die m\"{o}glichen Tracte auch als geschlossene ansehen,
und folglich aus jedem m\"{o}glichen durch Vorr\"{u}ckung des Anfangsgliedes andere ableiten.' (We can regard the
possible words also as closed curves, and therefore generate words by cyclic permutations of the labels.)

This results in 5 closed curves. He indicates the corresponding words in the table on the left side. 
The curves are drawn below this table (bottom of Figure \ref{Tractfiguren_handschriftlich}).
For instance, the words 1 and 39 (in the numbering from 1 through 105) belong to the first curve, and 54 and 89 to the second. 
For the third curve, associated with 4 words, Gauss sketched two variants. The remaining two curves each have 8 
representative words.

\begin{question}
What is the meaning of the table on the right? This contains the entries `I.4, IV.1, VIII.1, II.2, III.4, II.5' (in a triangular arrangement).
\end{question}

We now discuss the necessary condition for possible words and the simplified notation, which is similar to the Dowker-Thistlethwaite code.

\subsection{Even-odd rule and simplified notation}
We consider a possible word which gives a closed curve. Gauss noticed that the two occurrences of a label are separated by 
an even number of other labels. He explains that counting the occurrences yields one even and one odd number for each label. 
For the words numbered with 91 and 48 he gives examples, starting the counting with 0, and has the following correspondences:

$$ \text{For 91 (abcdbcda):} \left | \begin{array} {cc} 0 & 7 \\ 2 & 5 \\ 4 & 1 \\ 6 & 3 \end{array} \right |
\text{, and for 48 (abcabddc):} \left | \begin{array} {cc} 0 & 3 \\ 2 & 7 \\ 4 & 1 \\ 6 & 5 \end{array} \right | $$

\begin{figure}[hbtp]
\centering
\includegraphics[scale=0.44]{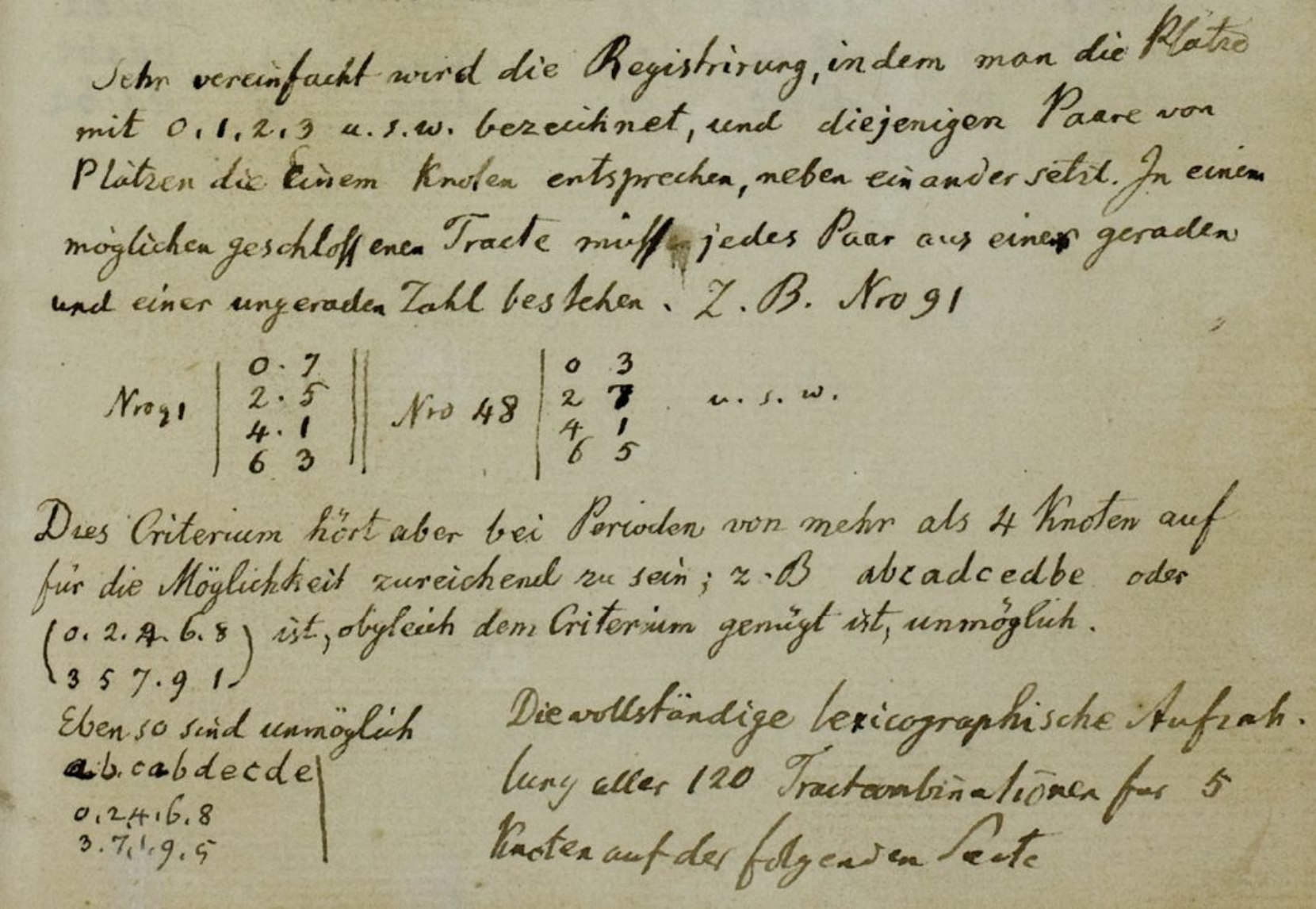}
\caption{The parity rule and the simplified notation} 
\label{Vereinfachte_Registrirung}
\end{figure}

The correspondence has a matrix form with ascending even numbers in the first and the respective odd numbers in the second column.
Tait's method is similar, and the Dowker-Thistlethwaite code for knot diagrams starts the counting with 1 and has odd numbers in 
the first column (or row if a transposed form is used).

We call the rule that labels are separated by an even number of other labels, the {\em Gauss parity rule}, 
and the matrix notation the {\em simplified notation}.

\clearpage
\newpage

Gauss explains (see Figure \ref{Vereinfachte_Registrirung}): 
`Sehr vereinfacht wird die Registrirung, indem man die Plätze mit $0, 1, 2, 3$ usw. bezeichnet, und
diejenigen Paare von Pl\"{a}tzen, die Einem Knoten entsprechen, neben einander setzt. In einem möglichen 
geschlossenen Tracte muss jedes Paar aus einer geraden und einer ungeraden Zahl bestehen. [\ldots] 
Dieses Criterium h\"{o}rt aber bei Perioden von mehr als 4 Knoten auf, f\"{u}r die M\"{o}glichkeit 
zureichend zu sein; z.B. abcadcedbe oder 
$\begin{psmallmatrix} 0 & 2 & 4 & 6 & 8\\ 3 & 5 & 7 & 9 & 1\end{psmallmatrix}$ 
ist, obgleich dem Criterium gen\"{u}gt ist, unm\"{o}glich. [\ldots]'.

\begin{figure}[hbtp]
\centering
\includegraphics[scale=0.43]{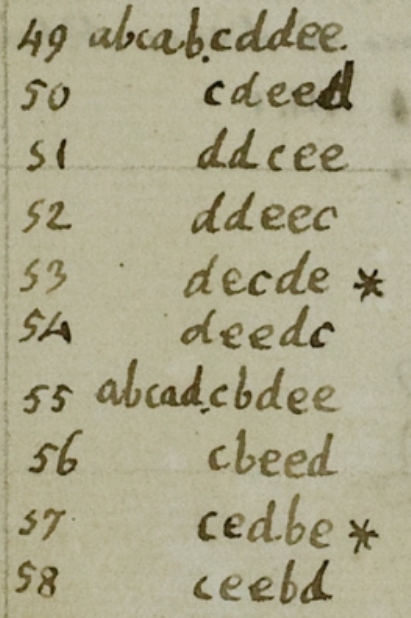} \quad
\includegraphics[scale=0.5]{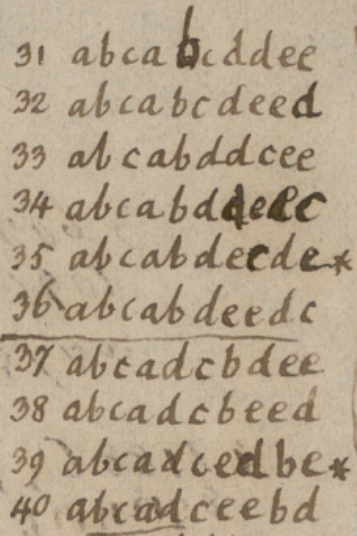} \quad
\includegraphics[scale=0.44]{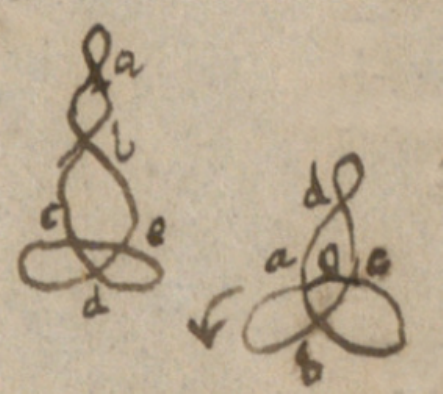}
\caption{Left: Gauss words for $n = 5$ in Handbuch 6, Middle: The same words in manuscript 33, Right: Sketches in manuscript 33} 
\label{gauss_n_5}
\end{figure}

\medskip
Translated into English this is:

\medskip
`The generation of possible words is much simplified when we number the places by $0, 1, 2, 3, \dots$
and write down the corresponding pairs. For a closed curve, belonging to a possible word, every pair
has to consist of an even and an odd number. [\ldots] For more than 4 crossings this criterion is no
longer sufficient; e.g. abcadcedbe or 
$\begin{psmallmatrix} 0 & 2 & 4 & 6 & 8\\ 3 & 5 & 7 & 9 & 1\end{psmallmatrix}$ 
is not possible, although it satisfies the condition. [\ldots]'

\subsection{Five crossings} 
Gauss gives another example for this phenomenon with abcabdecde, and he then writes, 
that the following page contains the 120 combinations for $n = 5$. He wrote down only those
words which satisfy his rule for possible words; and he certainly used the simplified notation
to generate the list (note that $5! = 120$). The list has only one type of star markings: The words
which cannot be realized as closed curves are marked with one star; there are 7 cases.
The two examples he gives in the text, and the star markings for them, are shown in Figure \ref{gauss_n_5} (left).

\subsection{A different manuscript page for the case of five crossings}
In the manuscript `Cod. Ms. Gauß Math. 33', \cite{Gauss_Ms_33}, consisting of only a few pages in an envelope, 
we find a different list of the 120 cases for $n = 5$ (see Figure \ref{gauss_n_5}, middle, for the same words as above).
The numbering is different (containing errors in the lexicografic ordering) but the markings with one star are also there.

The page also contains sketches of curves, see Figure \ref{gauss_n_5} (right).
The second curve is similar to the fifth curve in Figure \ref{Tractfiguren_handschriftlich}, 
but in Figure \ref{gauss_n_5} it has an additional curl and is marked with labels (the label for the curl is missing, however). 
To find such examples is interesting, because Gauss must have sketched  many cases in order to match permutations and words, 
and most of them are probably not conserved.

\clearpage
\newpage

\subsection{The contribution of Nagy}
We assume that Gauss wrote down these results around 1825; they were included in volume 8 of
his collected works, which appeared in 1900. The Hungarian mathematician Gyula Sz\"{o}kefalvi Nagy (1887--1953)
probably had access to this volume (he spent some time in G\"{o}ttingen in 1911/12) and published 
the first proof of the Gauss parity rule in 1925, and two years later in a German version, \cite{Nagy}.

\begin{figure}[hbtp]
\centering
\includegraphics[scale=0.6]{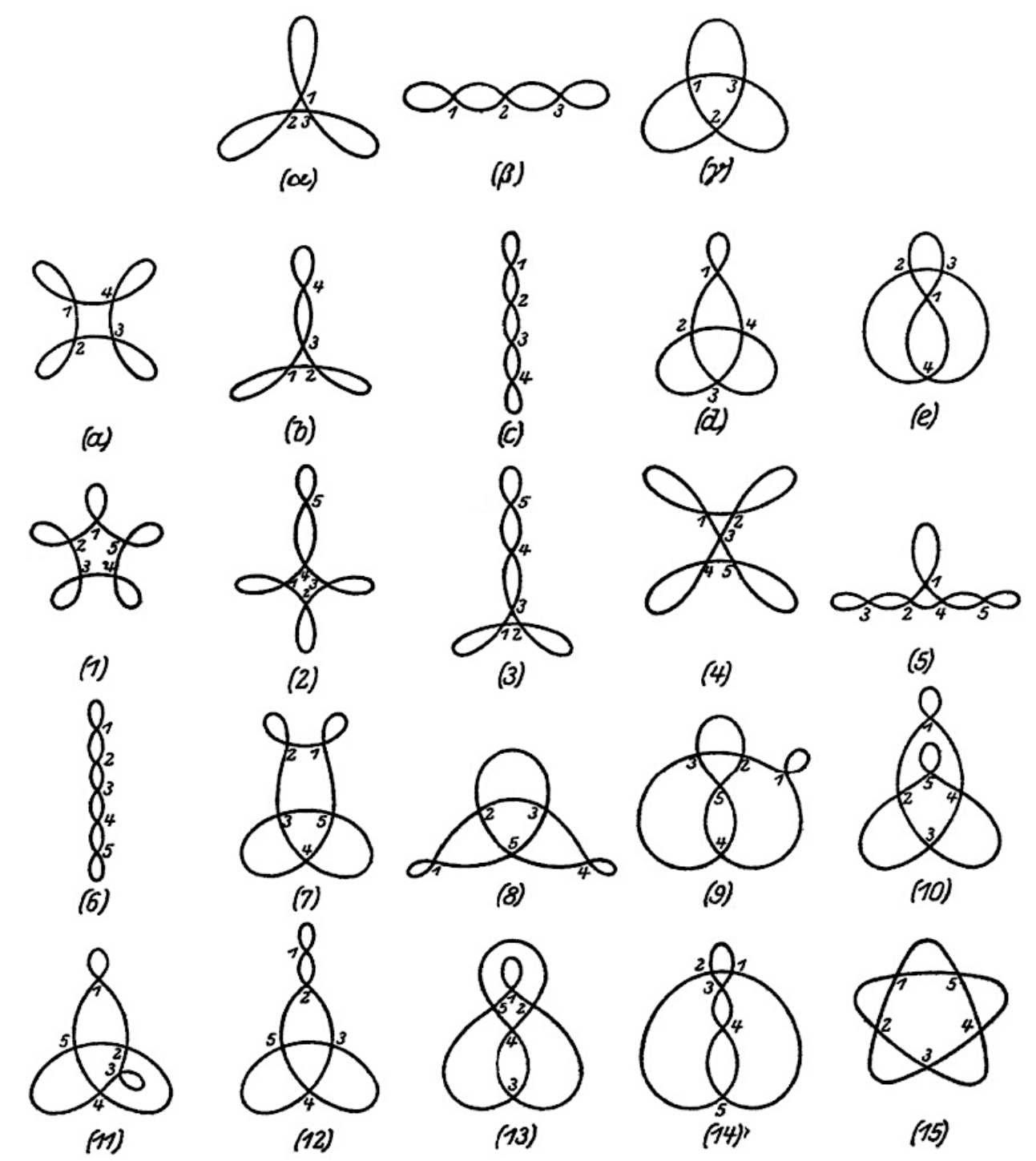}
\caption{Closed curves with $n = 3,4,5$ in Nagy \cite{Nagy}} 
\label{Nagy_Figure_1}
\end{figure}

The closed curves with $n = 3, 4, 5$ appear in his article, as shown in Figure \ref{Nagy_Figure_1}.
For $n=3$ his curves are labelled by $\alpha$, $\beta$, $\gamma$, for $n=4$ with $a,b,c,d,e$, and for $n=5$ with $1, \ldots, 15$.

Comparing the illustrations of closed curves of Gauss and Nagy, we find variations of the placement 
of curls, as we have shown in Figure \ref{template_n_2}. To see this, compare the curves (b) and (c) 
in Figure \ref{Nagy_Figure_1} with the curves drawn by Gauss. It is remarkable, that Gauss sketched the 
two variants of the third curve in Figure \ref{Tractfiguren_handschriftlich}, differing only in curl variations.

Shortly after, in 1936, Max Dehn found the first necessary and sufficient criterion for possible Gauss words 
(realizable as closed curves), see \cite{Dehn}.

\subsection{Arnold}
\begin{figure}[hbtp]
\centering
\includegraphics[scale=0.55]{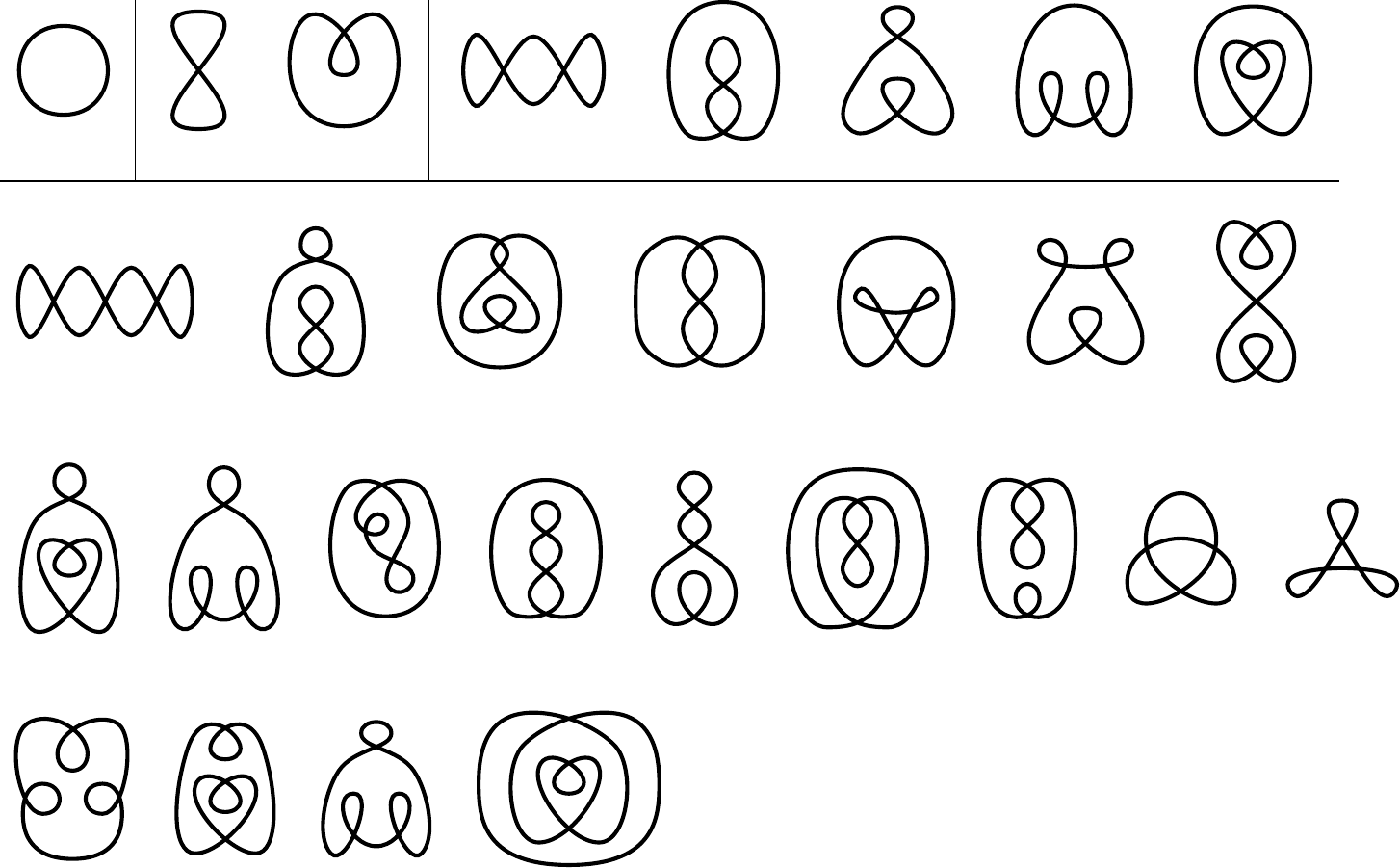}
\caption{Closed curves with $n \le 3$ in Arnold \cite{Arnold}} 
\label{Arnold}
\end{figure}

Vladimir Arnold studied closed and long curves in 1992, and he introduced the invariants $St$ (strangeness), $J^+$ and $J^-$. 
The curves in Figure \ref{Arnold} appeared in 1994 in \cite{Arnold}. We show his illustration for closed curves with $n \le 3$
because it allows a comparison with the figures given by Gauss and Nagy.

For the numbers of closed curves see also the sequence A008983 (Number of immersions of the unoriented circle into the unoriented plane with 
$n$ double points) in OEIS.
This sequence starts as follows:

\begin{center}
	\begin{tabular}[b]{|r|r|r|r|r|r|}
		0 & 1 & 2 & 3 & 4 & 5 \\
		\hline
		1 & 2 & 5 & 20 & 82 & 435 \\
	\end{tabular}
\end{center}

\medskip
Taking into account the orientation of the plane, we obtain the sequence A008981 (Number of immersions of the unoriented 
circle into the oriented plane with $n$ double points). It starts with $1, 2, 5, 21$. The different number for $n = 3$ is due
to the closed curve without vertical symmetry axis (find it!).

Arnold also studied long curves (in the next chapter in \cite{Arnold}).  
The number of long curves with $n$ crossings is given in the sequence A054993 in OEIS. 
For $n \le 4$ it has the values $1, 2, 8, 42, 260$. The 8 long curves for $n = 2$ are, 
for instance, depicted in Figure \ref{template_n_2} above. 
The closures of the 42 long curves for $n = 3$ lead to the 20 closed curves with $n = 3$ in Figure \ref{Arnold}.
The quite large number of 260 long curves for $n = 4$ indicates the difficulty we face with enumeration for 4 crossings.

\subsection{Comparison for the word abcddcba}
Nagy just lists one closed curve for the word abcddcba (this is (c) in Figure \ref{Nagy_Figure_1}).
Gauss was aware of different closed curves for the same word; he has two versions (third column with two curves in Figure \ref{Tractfiguren_handschriftlich}).
Arnold sketched all 26 closed curves (page 15 in \cite{Arnold}) and gives their invariants $St$, $J^+$ and $J^-$.

We conclude the historical digression with the remark, that although, as demonstrated above, Gauss developed his notation for closed and long curves, 
we did not yet find a sketch of a long curve in his hand-written documents.

\section{Conclusion of the case $n=2$} \label{result_n_2}
\begin{figure}[hbtp]
\centering
\includegraphics[scale=0.58]{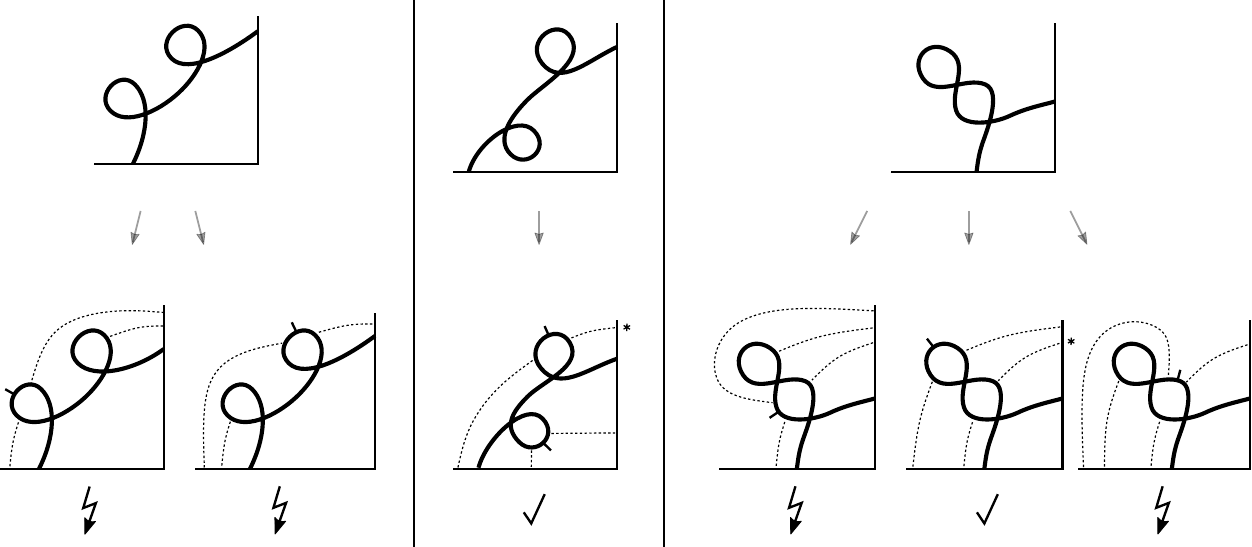}
\caption{The insertion of twist markers for long curves with $n = 2$, yielding 6 parterre templates, resulting in 2 valid templates after the filtering.}
\label{parterre_templates_n_2}
\end{figure}

We ended Subsection \ref{Case_n_two_firstPart} giving three long curves for further analysis.
These are shown in the top row of Figure \ref{parterre_templates_n_2}.
We describe the filter application to the parterre templates in the second row.
First, we introduce the concept of `separators'. These are markings on diagram arcs which separate 
twist markers going to the x-axis and those going to the y-axis.

\begin{definition} \label{separator}
In a parterre template, an \textit{outer separator} is a marking on an outer arc, so that all twist markers connected to an 
outer arc with lower index are twist markers to the x-axis, and all twist markers with larger index are twist markers to the y-axis,
see Figure \ref{template_regions} (index) and the lower row in Figure \ref{parterre_templates_n_2} (examples of separator markings). 
The marked arc may have both types of twist markers. An \textit{inner separator} is defined in the same way, but for inner arcs.
We allow only one outer and one inner separator in a parterre template.
\end{definition}

As a consequence of the `central arc flip' move M3, we do not connect twist markers to central arcs (see also principle P2).
Therefore, we do not need separators for central arcs.

We consider the first parterre template in the first row of Figure \ref{parterre_templates_n_2}.
There are two possibilities for outer separators, and the resulting parterre templates with twist markers are shown in the second row.
Note, that there is no possibility for a twist marker in the inner region, since in these templates all inner arcs are central arcs. 
The second template in the first row results in one parterre template (exactly one inner and one outer separator), and in the third case there are three possibilities.

The analysis of the cases follows:
\begin{enumerate}
	\item Filter F1 can be applied to the second curl. Therefore, this template would be reduced to a case with $n=1$, and we filter it out.
	\item Filter F4 (composite diagram parts) can be used at the first curl, and we disregard this case as well.
	\item This is a valid case. The twist marker with a star may be moved down (M2), however, and is ignored by M3. 
	      The parterre template without this twist marker is then used for further analysis as case X2 in Figure \ref{ParterreTemplatesNineCases}.
	\item Filter F1 can be applied to the outer curl. 
	\item A valid case, similar to the third case, resulting in X1 in Figure \ref{ParterreTemplatesNineCases}.
	\item This parterre template has only one twist marker to the y-axis. Filter F2 can be applied: It may be moved down (M2) and is ignored by M3.
	      The result is composite of the form $K \sharp -K$.
\end{enumerate}

We now consider the crossing variations for X1 and X2. Since these parterre templates each have two crossings, 
we would need to consider 4 variations in each case. However, by the `crossing choice' principle P3 we need only 2 variations in each case.
The resulting list of doubly symmetric diagrams has been analyzed with Knotscape, see the template results for $t_1$, $t_2$, $t_3$
in the Appendix in \cite {Lamm2023}. It turns out that for X1 only one variation is necessary, but for X2 both variations are needed
to generate minimal diagrams. The following theorem answers Question \ref{q2} and part of Question \ref{q1}.

\begin{theorem} \label{th14}
The prime knots with $c_{ds} = 14$ are $10_{99}$, $10_{123}$, \textnormal{12n706} and \textnormal{14n9732},
and there are no prime knots with $c_{ds} \le 12$.
\end{theorem}

\begin{proof}
We first note, that knot diagram templates with $n$ crossings lead to doubly symmetric diagrams with at least $4n + 4$ crossings, 
if the resulting knots are required to be prime. The reason is, that if the x- or y-axis does not contain a crossing, the diagram 
is composite. Hence, for non-trivial prime knots there must be at least 2 crossings on each of the axes. If $n \ge 3$, the diagrams 
therefore have at least 16 crossings, and the analysis of templates X1 and X2 is sufficient for $c_{ds} \le 14$. 
Using Knotscape, results in the statement of the theorem.
\end{proof}

\begin{figure}[hbtp]
\centering
\includegraphics[scale=0.8]{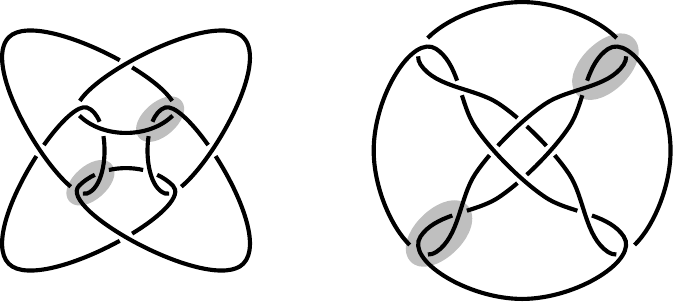}
\caption{Left: The knot $t_1(1,0 \mid 1) = 4_1 \sharp 4_1$ with partial knot $5_1$.
This doubly symmetric diagram with 12 crossings represents a composite knot.
Right: The same diagram, but with switched inner and outer regions.}
\label{crossings_12_composite}
\end{figure}

Figure \ref{crossings_12_composite} shows on the left side a diagram with 12 crossings generated with template X1.
It has 4 crossings on the axes and we see a cancelling of twists on a diagonal (marked in gray).
Although this diagram does not obviously show a composite knot, it is composite.
Diagrams with up to 12 crossings seem to lack the complexity which can prevent such twist cancellations.
On the right side, we switched inner and outer regions and the twist cancelling is seen in a more obvious way there.

\section{The generation and filtering of templates for $n=3$} \label{result_n_3}
\begin{figure}[ht]
\centering
\includegraphics[scale=0.55]{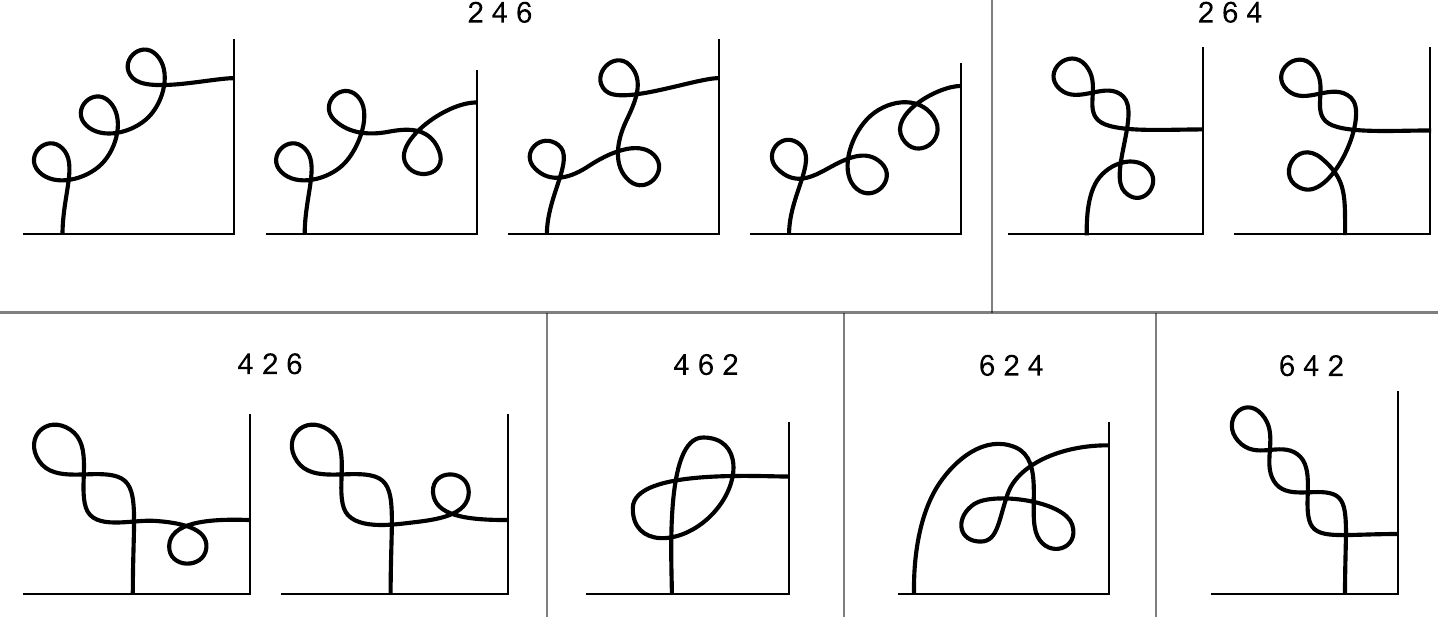}
\caption{For the 6 permutations with $n=3$ we generate all long curves for our enumeration goal.}
\label{parterre_templates_n_3}
\end{figure}

We mentioned that for $n=3$ there are 42 long curves. If we disregard those curves with Reidemeister 1 curls in the interior of a 
diagram region (filter F5) and choose only one curve from each pair related by the switching of inner and outer regions, then 11 curves remain;
a possible choice is shown in Figure \ref{parterre_templates_n_3}. We have included the Dowker-Thistlethwaite codes (the six permutations
of the even numbers 2, 4, 6), which are more economical than Gauss words---as Gauss already noted.

\begin{figure}[hbtp]
\centering
\includegraphics[scale=0.6]{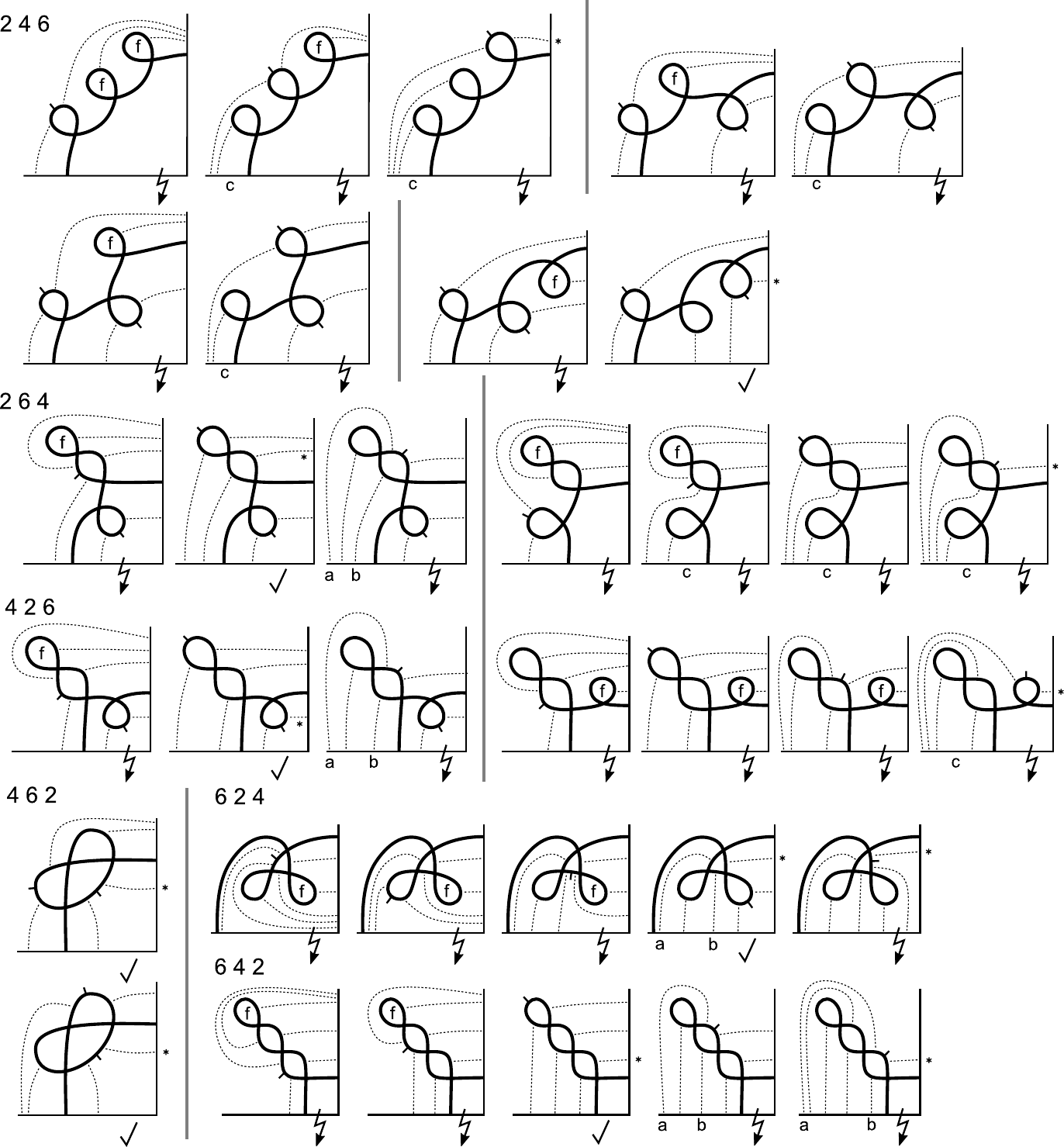}
\caption{The 35 templates with twist markers for $n=3$. After filter applications 7 valid templates remain.}
\label{parterre_templates_n_3_with_twist_markers}
\end{figure}

Proceeding in the same way as for $n=2$, we list all possibilities for the outer and inner separators and analyze
the twist marker configurations, see Figure \ref{parterre_templates_n_3_with_twist_markers}.

\clearpage
\newpage

We explain the application of the filters, and the markings in Figure \ref{parterre_templates_n_3_with_twist_markers} with 
`f', `*', `c' and `a b':

\begin{itemize}
	\item The curls, for with the F1 filter is applied, are marked with `f' inside the curl (17 templates are filtered out).
  \item Applied to the remaining 18 templates: If a template contains only one twist marker connected to the y-axis, and this can be 
	      removed by filter F2, we mark it with a star (5 templates are filtered out).
	\item With filter F4 at the regions marked with `c', we remove 3 additional templates.
	\item From the remaining 10 templates, 3 can be filtered out with filter F3 (with markings `a' and `b').
	      In these 3 cases, one twist marker can be omitted. We then conclude that each case is reduced to one of the 7 remaining cases.
\end{itemize}

Each of the 7 valid templates contains a star marker. We explain its meaning:
In the two cases for permutation 462 the twist markers with star can be moved up (M2 move) and then added to the 
twist marker on the other side of the curve's end point at the y-axis (as illustrated in Figure \ref{filter_examples}). 
In the other 5 cases the twist markers with star can be moved up or down and then omitted, because they end at a central arc (filter F2). 
After these simplifications we arrive at the cases with $n=3$ (H1, $\ldots$, H7) in Figure \ref{ParterreTemplatesNineCases}.

We consider the crossing variations for H1, $\ldots$, H7; we need 4 variations in each case.
The resulting list of doubly symmetric diagrams has been analyzed with Knotscape, see the template results for $t_4, \ldots t_{16}$
in the Appendix in \cite {Lamm2023}. This list contains all templates which are necessary to generate minimal diagrams
(note, for instance, that there are no templates corresponding to $H_1$).
The following theorem answers the remaining part of Question \ref{q1}.

\begin{theorem}
There are 17 prime knots with $c_{ds} = 16$ and 26 prime knots with $c_{ds} = 18$. These are given in the following table.
\end{theorem}

\small 
\begin{center}
	\begin{tabular}[b]{l|lllllll}
		 $c_{ds}$ &&&&&&& \\
		 \hline
		 16 & $12a_{427}$ & $12a_{1019}$ & $12a_{1105}$ & $12a_{1202}$ &&& \\
		    & $14a_{8662}$ & $14a_{18676}$ & $14a_{19472}$ &&&& \\
				& $14n_{8213}$ & $14n_{22073}$ & $14n_{25903}$ &&&& \\
				& $16n_{428839}$ & $16n_{451788}$ & $16n_{645918}$ & $16n_{645926}$ & $16n_{847920}$ & $16n_{991381}$ & $16n_{991505}$ \\
		 \hline
		 18 & $14a_{6002}$ & $14a_{16311}$ & $14a_{17173}$ & $14a_{18187}$ & $14a_{18362}$ & $14a_{18680}$ & $14a_{18723}$ \\
		    & $16a_{107430}$ & $16a_{313024}$ & $16a_{314171}$ & $16a_{330218}$ & $16a_{354511}$ & $16a_{356843}$ & \\
        & $16n_{101996}$ & $16n_{102000}$ & $16n_{102453}$ & $16n_{106013}$ & $16n_{268599}$ & $16n_{323632}$ & $16n_{797553}$ \\
				& $16n_{847983}$ & $16n_{858257}$ & $16n_{868471}$ & $16n_{869383}$ & $16n_{872167}$ & $16n_{872172}$ & \\
		 \hline
	\end{tabular}
\end{center}
\normalsize

\begin{proof} 
Analogously to the proof of Theorem \ref{th14}, if $n \ge 4$, the diagrams have at least 20 crossings, and the analysis of the templates X1, X2, H1, $\ldots$, 
H7 is sufficient for $c_{ds} = 16$ and $c_{ds} = 18$. Using Knotscape for the crossing variations results in the statement of the theorem.
\end{proof}

\section{The generation and filtering of templates for $n=4$}
For $n=4$ we proceed in the same way as for $n=3$, but encounter a lot more cases:
There are $24$ Dowker-Thistlethwaite codes and 260 long curves. Not taking into account curves with Reidemeister 1 curls in the 
interior of a diagram region, and choosing only one curve from each pair related by the switching of inner and outer regions (P1),
yields 49 long curves. 

\begin{figure}[hbtp]
\centering
\includegraphics[scale=0.65]{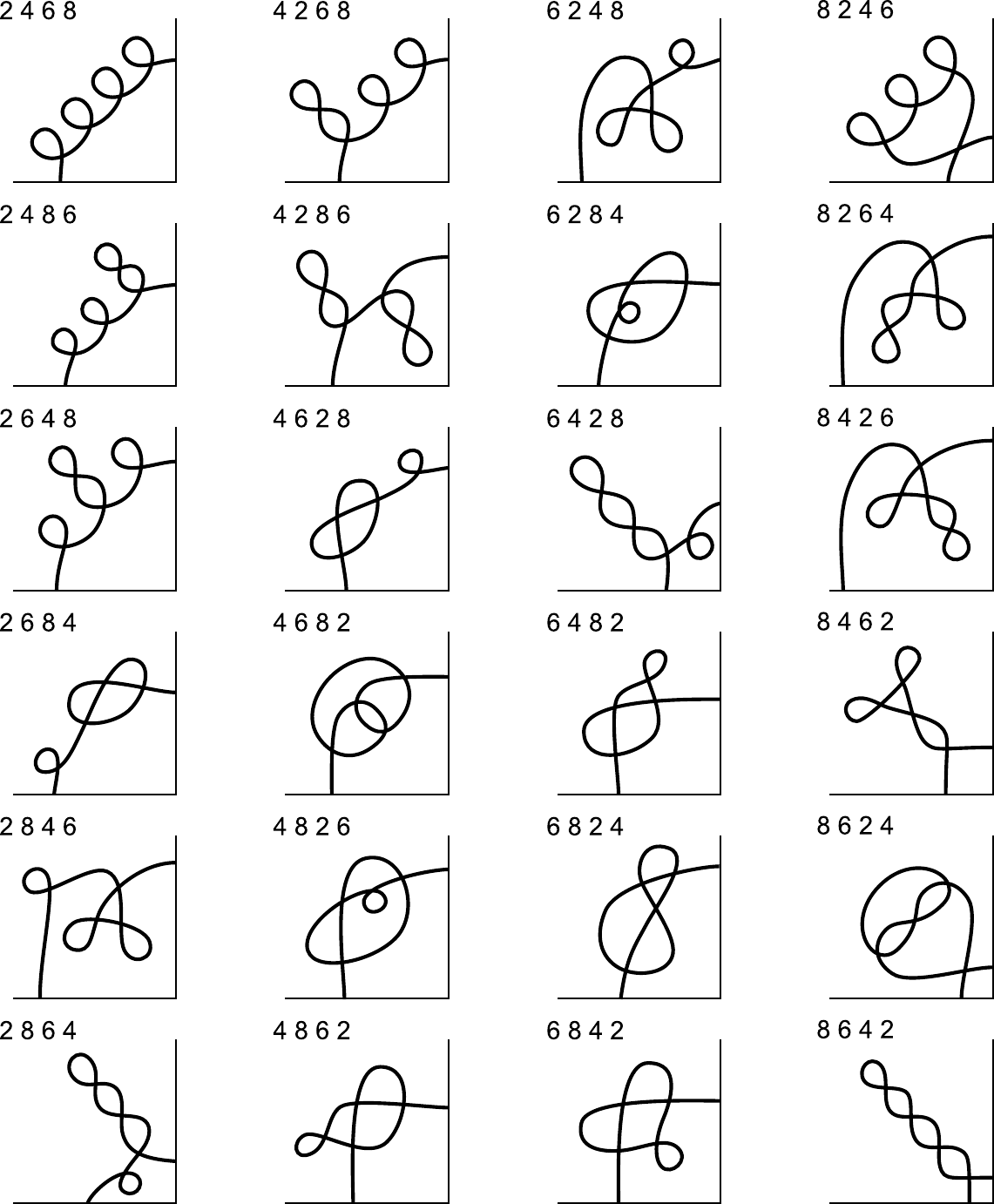}
\caption{The Dowker-Thistlethwaite codes for $n=4$, showing one example template in each case.}
\label{parterre_templates_n_4}
\end{figure}

Figure \ref{parterre_templates_n_4} contains one example template for each permutation of the even numbers 2, 4, 6, 8.
In Table \ref{curve_counts_n_4} we summarize the long curve properties.

Column `DT' contains the Dowker-Thistlethwaite code, column $m_1$ the number of long curves with this code. In the
next column, with header $m_2$, we give the number of long curves without Reidemeister 1 curls in the interior of a diagram region (filter F5).
The last column, $m_3$, is the result of choosing only one curve from each pair related by the switching of inner and outer regions (P1),
and we have $m_3 = m_2 / 2$. The sum of all entries in the $m_1$ columns is 260, the sums of the entries in the following columns are 98, and 49.

Compared to the 11 cases for $n=3$ in Figure \ref{parterre_templates_n_3}, the number of cases for $n=4$ has increased considerably.
Therefore, for the filtering of the cases including the placement of separators, for $n=4$ we list only the remaining cases,
after the filters have been applied. The 37 remaining cases are shown in Figure \ref{parterre_templates_n_4_with_twist_markers}. 

\small
\begin{table}[hbtp]
	\centering
		\begin{tabular}[b]{|r|r|r|r||r|r|r|r||r|r|r|r||r|r|r|r|}
			\hline
			DT & $m_1$ & $m_2$ & $m_3$ & DT & $m_1$ & $m_2$ & $m_3$ & DT & $m_1$ & $m_2$ & $m_3$ & DT & $m_1$ & $m_2$ & $m_3$\\ \hline
			2468 &  16 & 16 & 8 & 4268 & 16 & 8 & 4 & 6248 &  16 &  4 & 2 & 8246 & 16 & 2 & 1 \\ \hline
			2486 &  16 &  8 & 4 & 4286 & 16 & 4 & 2 & 6284 &   4 &  0 & 0 & 8264 & 16 & 2 & 1 \\ \hline
			2648 &  16 &  8 & 4 & 4628 &  4 & 4 & 2 & 6428 &  16 &  4 & 2 & 8426 & 16 & 2 & 1 \\ \hline
			2684 &   4 &  4 & 2 & 4682 &  2 & 2 & 1 & 6482 &   4 &  4 & 2 & 8462 & 16 & 2 & 1 \\ \hline
			2846 &  16 &  4 & 2 & 4826 &  4 & 0 & 0 & 6824 &   2 &  2 & 1 & 8624 &  4 & 4 & 2 \\ \hline
			2864 &  16 &  4 & 2 & 4862 &  4 & 4 & 2 & 6842 &   4 &  4 & 2 & 8642 & 16 & 2 & 1 \\ \hline
		\end{tabular}
	\caption{Properties of the long curves for $n=4$.}
	\label{curve_counts_n_4}
\end{table}

\normalsize

\subsection{Irreducible curves have $m_1=2$}
Arnold defines reducible closed curves (p.\,18 in \cite{Arnold}). We follow this definition and apply it to long curves, using their closure. 
The closure of a long curve (in our template form) is obtained by connecting the start and end points along the x- and y-axis.
A long curve is called \textit{reducible} if some crossing cuts the closure of the long curve into two disjoint parts.
If a long curve is not reducible, we call it irreducible.

Out of the 24 long curves in Figure \ref{parterre_templates_n_4} there are only two irreducible curves: 4682 and 6824.
The other 22 long curves can be generated as the addition of two curves, or by the more general operation of inserting
a long curve at some arc of another curve. The simplest case for this is the insertion of a Reidemeister 1 curl, as seen
in the cases 4826 and 6284. Because these two cases each have a curl in the interior of a diagram region, they have $m_2=0$. 

For an irreducible long curve, the number $m_1$ is equal to 2. If the curve splits into two irreducible parts, 
then $m_1=4$, and otherwise $m_1=16$. We observe, that for $n=4$ the case of three irreducible parts does not occur, 
because there is no irreducible long curve with two crossings.

\subsection{Filter application}
\begin{figure}[hbtp]
\centering
\includegraphics[scale=0.57]{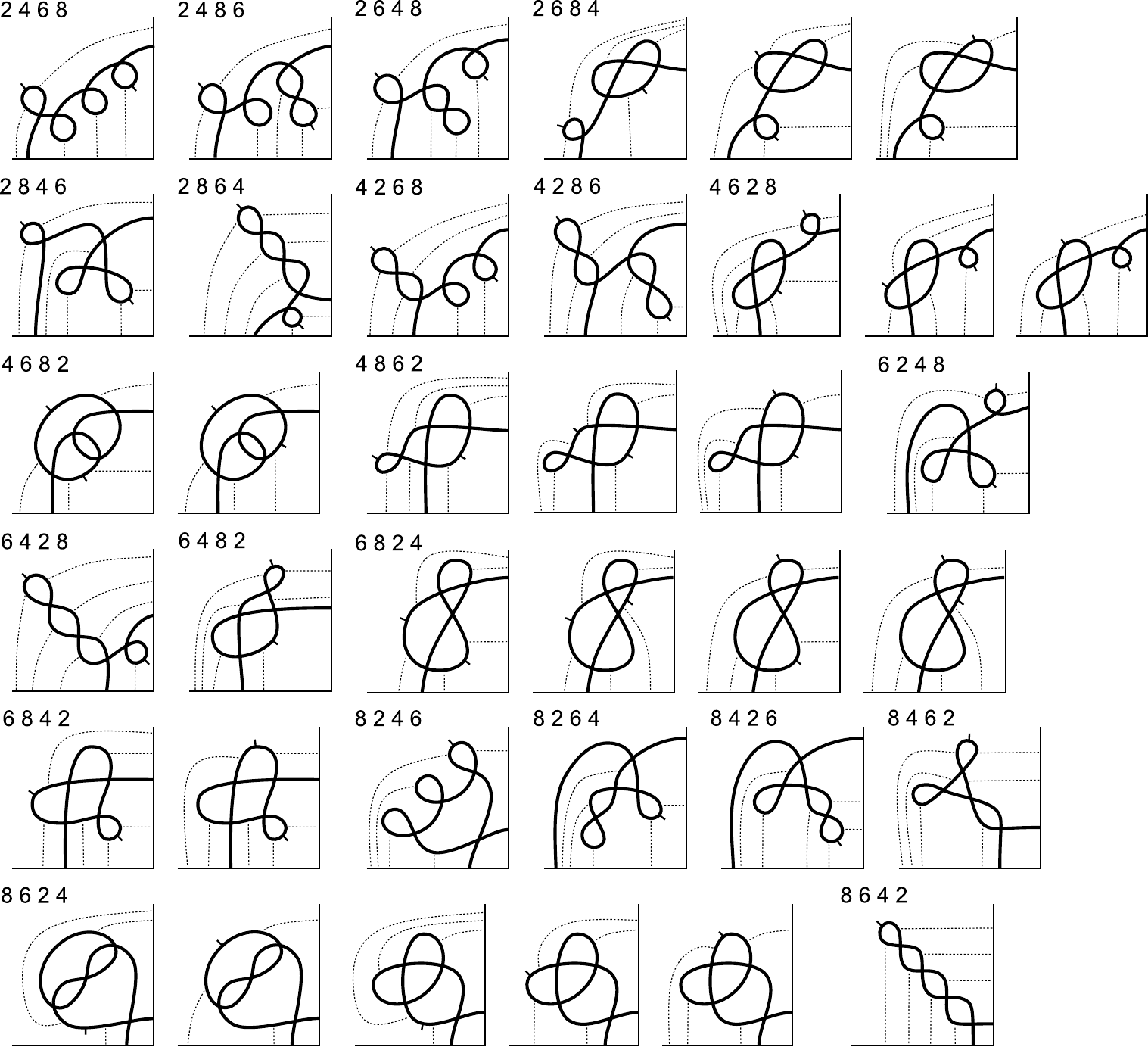}
\caption{The 37 templates with twist markers for $n=4$ (after filters have been applied).}
\label{parterre_templates_n_4_with_twist_markers}
\end{figure}

The application of filters is similar to what we described in more detail for $n \le 3$ and we leave the check, that 
this results in the 37 cases in Figure \ref{parterre_templates_n_4_with_twist_markers}, as an exercise.

Again, we consider the crossing variations for these templates; we need 8 variations in each case.
Due to the manual effort of encoding the templates as program input, we used only 15 of the templates for further analysis.
The resulting list of doubly symmetric diagrams has been analyzed with Knotscape. 
Minimal diagrams are generated by the templates $t_{17}, \ldots t_{23}$ in the Appendix in \cite {Lamm2023}. 

\subsection{Properties of templates with crossing variations}
\begin{figure}[hbtp]
\centering
\includegraphics[scale=0.8]{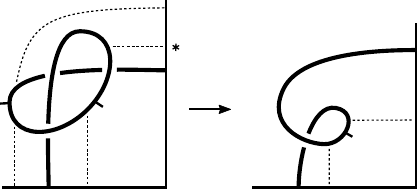}
\caption{Filtering through simplification can occur also for templates with crossing variations.}
\label{crossing_variations}
\end{figure}

Simplifications can also occur for templates with crossing variations. If this is the case, we omit the knot diagram template
from further analysis. For an example, see Figure \ref{crossing_variations}, where a Reidemeister 2 move is possible. 
As a result, we do not need all of the $8 \times 37$ crossing variations (and we also utilized only a subset of the 
$8 \times 15$ templates which were in fact taken into account).

\section{Summary and outlook}
Our template approach has been successful for $n \le 3$. Starting with $n=4$, the number of possible templates is too large
to be processed in an analysis which includes manual steps. Therefore, a different approach, or an automated processing
would be necessary. On the other hand, even a complete generation of all doubly symmetric diagrams up to a certain number of
crossings (e.g.\;26 or 28) cannot solve the question, which of the remaining (`almost doubly symmetric') ribbon knots up to
16 crossings are doubly symmetric. 

The following description, taken from the first part of the article, is valid:
The main open problem is to find a knot invariant which enables us to show that there is a strongly positive 
amphicheiral knot which is ribbon, or even a symmetric union, but does not have a doubly symmetric diagram.

%\clearpage
%\newpage

\section*{Acknowledgments}
I thank Michael Eisermann for his help in the interpretation of Gauss's star markings.

%%%%%%%%%%%%%%%%%%%%%%%%%%%%%%%%%%%%%%%%%%%%%%%%%%%%%%%

\vspace{1cm}
\noindent
Christoph Lamm \\ \noindent
R\"{u}ckertstr. 3, 65187 Wiesbaden \\ \noindent
Germany \\ \noindent
e-mail: christoph.lamm@web.de

\end{document}